\documentclass[11pt,titlepage,a4paper]{article}
\usepackage{bozhomac}	
\usepackage{bozhlogo}	
\usepackage{cite}	

\title{\bfseries	\vspace*{-1.678902345in}
{\huge   Deviation equations\\[1ex]
	 in spaces with affine connection        }
}

\vspace{1.7ex}

\author{
Bozhidar Z.\ Iliev
\thanks{Laboratory of Mathematical Modeling in Physics,
Institute for Nuclear Research and \mbox{Nuclear} Energy,
Bulgarian Academy of Sciences,
Boul.\ Tzarigradsko chauss\'ee~72, 1784 Sofia, Bulgaria}
\thanks{E-mail address: bozho@inrne.bas.bg}
\thanks{URL: http://theo.inrne.bas.bg/$\sim$bozho/}
\\
\fbox{Sawa S.\ Manoff}
\thanks{Laboratory of Solitons, coherence and geometry,
Institute for Nuclear Research and \mbox{Nuclear} Energy,
Bulgarian Academy of Sciences,
Boul.\ Tzarigradsko chauss\'ee~72, 1784 Sofia, Bulgaria}
\thanks{Prof.\ Sawa Manoff passed away on May 27, 2005}
\thanks{URL: http://theo.inrne.bas.bg/elpart/SavaManov/smanoff.html}
}

\date{
 \vspace{2.27ex}\ShortTitle{Deviation equations in spaces with
				affine connection}	\\[0.27ex]
 \vspace{3.27ex}
\small
	\begin{tabular}{r@{$\colon\to~$}l}
 \vspace{0.09ex} Basic ideas	& 1980--1983	\\[0.09ex]
\vspace{0.09ex} Initial typeset& 1983 	\\[0.09ex]
\vspace{0.09ex} Secondary typeset& November 2--23, 2005	\\[0.09ex]
 \vspace{0.09ex} Last update	& December 14, 2005	\\[0.09ex]
 \vspace{0.27ex} Produced	& \fbox{\today}	\\[0.27ex]
	\end{tabular} \\[1.27ex]
\normalsize
	\begin{tabular}{r@{$\colon~$}l}
 \vspace{0.27ex}
 \normalsize\sffamily\bfseries
  http://arXiv.org e-Print archive No. &
 \normalsize\sffamily\bfseries  math.DG/0512008 \\[1.27ex]
\vspace{0.27ex}
\bfseries Published as \bfseries\itshape
JINR Communication  P2-83-897 &
\bfseries\itshape JINR, Dubna, 1983 (In Russian)
   							 \\[0.27ex]
\vspace{0.27ex}
\bfseries Translation from Russian and new typeset &
\bfseries\itshape Bozhidar Z. Iliev
   							 \\[0.27ex]
	\end{tabular} \\[-0.27ex]
 \vspace{4.27ex}{\Huge\BOZHO}	\\[4.27ex]
 \vspace{0.27ex}\Subject{General relativity, Differential geometry} \\[2.27ex]
	\begin{tabular}{r@{\hspace{0.512em}}|@{\hspace{0.512em}}l}
 \vspace{0.27ex}	
&
 \vspace{0.27ex} 
	\end{tabular} \\[1.27ex]
\vspace{0.27ex}
 							\\[0.27ex]
}

\begin{document}		

\renewcommand{\thepage}{\roman{page}}

\renewcommand{\thefootnote}{\fnsymbol{footnote}} 
\maketitle				
\renewcommand{\thefootnote}{\arabic{footnote}}   



 	\begin{abstract}
Connections between Lie derivatives and the deviation equation has been
investigated in spaces $L_n$ with affine connection. The deviation equations
of the geodesics as well as deviation equations of non-geodesics trajectories
have been obtained on this base. This is done via imposing certain conditions
on the Lie derivatives with respect to the tangential vector of the basic
trajectory.
 	\end{abstract}

\renewcommand{\thepage}{\arabic{page}}



\addtocounter{section}{-1}

\section {Introduction}
\label{Introduction}

	In the last years, the deviation equations find a broad application
in the study of many mathematical and physical problems in gravitational
physics and astrophysics~\cite{Manov-1978,MTW,Mitskevich}. Several versions
of these equations are known in Riemannian spaces $V_n$.
In~\cite{Manov-1978}, on the base of Lie derivatives, a general method is
proposed for derivation of deviation equations in $V_n$.

	The purpose of the present investigation is in the application of the
methods in~\cite{Manov-1978} to spaces (manifolds) with affine connection
$L_n$. The approach proposed for finding deviation equations in  $L_n$ is
illustrated on a number of examples.

	The so-called \emph{generalized deviation equation} in spaces with
affine connection is considered in Sect.~\ref{Sect1}. In Sect.~\ref{Sect2}, a
physical interpretation of that equation is presented. Sect.~\ref{Sect3}
contains examples of additional conditions restricting the form of the
generalized deviation equation. Particular examples of deviation equations
are given in Sect.~\ref{Sect4}. In the Appendix, Sect.~\ref{Appendix}, is
recalled  the notion of Lie derivative of the connection coefficients of a
linear connection~\cite{Yano/LieDerivatives}.


\section {The generalize deviation equation in\\
					spaces with affine connection}
\label{Sect1}

\textbf{1.1.}
	The considerations that follow below are based on the identity
(see the Appendis, section~\ref{Appendix})
	\begin{equation}	\label{1.1}
\Lied_\xi\Gamma_{ij}^{k}
=
\xi_{|i|j}^{k} - R_{ijl}^{k} \xi^l - (T_{il}^{k}\xi^l)_{|j}.
	\end{equation}
Here all quantities are evaluated at a point $x\in L_n$ and the following
notation is introduced:
 $\Lied_\xi$ is the Lie derivative operator~\cite{Yano/LieDerivatives} along a
$C^2$ vector field $\xi=\xi^i E_i$;
the Latin indices $i,j,k,\dots$ run from 1 to $n\in\field[N]$ and summation
 over indices repeated on different levels is understood over the whole range
 of their values;
 $\{E_i\}$ is a frame in the tangent bundle $T(L_n)$, \ie
  $E_i=A_{i}^{\alpha}\pd_\alpha$, where the Greek indices
 $\alpha,\beta,\dots =1,\dots,n$ number a coordinate frame
 $\bigl\{\pd_\alpha:=\frac{\pd}{\pd x^\alpha}\bigr\}$
 and $A=[A_{i}^{\alpha}]$ is a non-degenerate matrix-valued function;
  $\Gamma_{ij}^{k}$ are the coefficients of the affine connection of $L_n$ in
the frame $\{E_i\}$ and are generally non\ndash symmetric in their
subscripts;
$\xi^k_{|i|j}:=\xi^k_{|ij}:=(\xi^k_{|i})_{|j}$ with the suffix ``$|j$'' denoting
the action of the covariant derivation operator along the basic vector field
$E_j$ relative to the connection with coefficients $\Gamma_{ij}^{k}$;
respectively with the suffix ``$,i$'' it will be denoted the action of $E_i$ on
the components of objects over $L_n$ considered as scalar functions, for
instance, we have $\xi_{|i}=\xi_{,i}+\Gamma_{kj}^{i}\xi^k$ with
$\xi^k_{,i}:=E_i(\xi^k)$;
\[
R_{jkl}^{i}
:=
- 2\Gamma_{j[k,l]}^{i} - 2_{j[k}^{n} \Gamma_{|n|l]}^{i} -
\Gamma_{jn}^{i} C_{kl}^{n}
\]
are the components of the curvature tensor in $\{E_i\}$, where
	\begin{gather*}
B_{[ij]} ;= \frac{1}{2}(B_{ij}-B_{ji})\qquad
B_{[i|j|l]} ;= \frac{1}{2}(B_{ijl}-B_{lji})
\\
C_{jk}^{i} := -2 A_{\alpha}^{i} A_{[j,k]}^{\alpha}\qquad
[A_{\alpha}^{i}] := [A_{j}^{\beta}]^{-1};
	\end{gather*}
and
\[
T_{jk}^{i} := -2 \Gamma_{[jk]}^{i} - C_{jk}^{i}
\]
are the components of the torsion tensor of $L_n$ relative to $\{E_i\}$,

\textbf{1.2.}
	Let $u=u^iE_i$ be arbitrary contravariant  $C^1$ vector field.
Projecting~\eref{1.1} on $u^i$  and $u^j$ and taking into account the
equation
$ u^iu^j \xi_{|ij} = u^j(\xi^k_{|i}u^i)_{|j} - \xi^k_{|i}u_{|j}^{i} u^j $,
we get the relation
	\begin{equation}	\label{1.2}
u^j (u^i\xi_{|i}^k)_{|j}
=
R_{ijl}^{k}u^iu^j\xi^l +
\xi_{|i}^{k} (u_{|j}^{i}u^j) +
u^iu^j (T_{il}^{k}\xi^l)_{|j} +
u^iu^j\Lied_\xi \Gamma_{ij}^{k} .
	\end{equation}

	If $u$ is a $C^2$ vector field, then
\[
u^iu^j\Lied_\xi \Gamma_{ij}^{k}
=
\Lied_\xi(u^k_{|i} u^i) - u^i(\Lied_\xi u^k)_{|i} - u_{|i}^{k}\Lied u^i
\]
and, consequently,~\eref{1.2} can be rewritten as
	\begin{equation}
\tag{\ref{1.2}a}	\label{1.2a}
u^j (u^i\xi_{|i}^k)_{|j}
=
R_{ijl}^{k}u^iu^j\xi^l +
\xi_{|i}^{k} (u_{|j}^{i}u^j) +
u^iu^j (T_{il}^{k}\xi^l)_{|j} +
\Lied_\xi(u^k_{|i} u^i) - u^i(\Lied_\xi u^k)_{|i} - u_{|i}^{k}\Lied u^i .
	\end{equation}

\textbf{1.3.}
	Let $x(s)=(x^1(s),\dots,x^n(s))$ be a $C^2$ path with parameter $s$
in $L_n$ and the vector field $u$ be chosen such that
	\begin{equation}	\label{1.3}
u^\alpha\big|_{x=x(s)}
=  \frac{\od x^\alpha(s)}{\od s}
:= u^\alpha(s) =u^\alpha ,
	\end{equation}
\ie $u(s)$ to be the vector tangent to $x(s)$ at the parameter values $s$.

	Denote by $\bar{\nabla}_i$ the covariant derivative operator of
tensor field components; for example $\bar{\nabla}_i\xi^k=\xi_{|i}^{k}$. Then
\(
\frac{\bar{D}}{\od s}
:=u^i\bar{\nabla}_i
=u^\alpha\bar{\nabla}_\alpha
=\frac{\od x^\alpha(s)}{\od s}\bar{\nabla}_\alpha
\)
is the covariant derivative along $x(s)$ with respect to $s$.

	Using the operator $\bar{D}/\od s$, we can rewrite~\eref{1.2}
and~\eref{1.2a} along $x(s)$ respectively as
	\begin{align}	\label{1.4}
\frac{\bar{D}^2\xi^k}{\od s^2}
& =
R_{ijl}^{k} u^iu^j\xi^l +
\xi_{|j}^{k} F^j +
u^j \frac{ \bar{D} (T_{jl}^{k}\xi^l) } {\od s} +
u^iu^j \Lied_\xi \Gamma_{ij}^{k},
\\
\tag{\ref{1.4}a}	\label{1.4a}
\frac{\bar{D}^2\xi^k}{\od s^2}
& =
R_{ijl}^{k} u^iu^j\xi^l +
\xi_{|j}^{k} F^j +
u^j \frac{ \bar{D} (T_{jl}^{k}\xi^l) } {\od s} +
\Lied_\xi F^k - \frac{\bar{D}(\Lied_\xi u^k)}{\od s} - u^k_{|i} \Lied_\xi u^i,
	\end{align}
where
\[
\frac{\bar{D}^2}{\od s^2}
:=
\frac{\bar{D}}{\od s} \frac{\bar{D}}{\od s},
\qquad
F^k:=\frac{\bar{D}u^k}{\od s} = u_{|i}^{k}u^i .
\]

	The equation~\eref{1.4} or its equivalent version~\eref{1.4a} will be
called \emph{generalized deviation equation} in spaces with affine connection
$L_n$ in the frame $\{E_i\}$. By its essence equation~\eref{1.4} is an
identity, which is a corollary of the definitions of the quantities and
operators entering in it. It can be regarded as an equation relative to some of
the quantities in it only if it is presupposed that between them some
connection exists. So, rigorously speaking, by a deviation equation we shall
understand the identity~\eref{1.4} together with some additional condition(s)
imposed on the quantities entering in it. Namely the variety of the possible
additional conditions is the cause for the existence of different forms
(modifications) of deviation equations; see~\cite{Manov-1978} in a case of
a Riemannian manifold $V_n$. It should be mentioned, the additional conditions
are sometimes useful to be considered as first integrals of the deviation
equation they determine.

	\begin{Rem}	\label{Rem2.1}
In a coordinate/frame-independent language, the equations~\eref{1.4}
and~\eref{1.4a} respectively read:
	\begin{align}	\label{1.5}
\frac{D^2\xi}{\od s^2}
&=
\Hat{R}(u,\xi)u + C_{1}^{1}(F\otimes D\xi) +
\frac{D(\Hat{T}(F,\xi))}{\od s} - \Hat{T}(F,\xi) +
\Lied_\xi\Gamma(u,u)
\\
\tag{\ref{1.5}a}	\label{1.5a}
	\begin{split}
\frac{D^2\xi}{\od s^2}
& =
\Hat{R}(u,\xi)u + C_{1}^{1}(F\otimes D\xi) +
\frac{D(\Hat{T}(F,\xi))}{\od s} - \Hat{T}(F,\xi) +
\Lied_\xi F -
\\ \hphantom{ \frac{D^2\xi}{\od s^2}= =} &\quad
-\frac{D(\Lied_\xi u)}{\od s} -
C_{1}^{2}\{(Du)\otimes (\Lied\xi u))\} .
	\end{split}
	\end{align}
Here we have introduce the following notation:
$\frac{D}{\od s} := u^i\nabla_i$ with $\nabla_i=\nabla_{E_i}$ is the
covariant derivative along $E_i$ (\eg $\nabla_i u=(\bar{\nabla}_i u^k)E_k$);
 $\xi=\xi^kE_k$,  $U=u^kE_k$,  $F=F^kE_k$;
$\Hat{R}(u,\xi)=[\nabla_u,\nabla_\xi] -\nabla_{[u,\xi]} $,
with $\nabla_u=u^i\nabla_i$, is the curvature operator in $L_n$;
 $[A,B]=AB-BA$ is the commutator of $A$ and $B$;
 $C_{q}^{p}$ is the contraction operator over the $p$\ndash superscript and
$q$\ndash subscript (recall $[C_{q}^{p},\nabla]=0$);
 $DA:=(\nabla_k A)\otimes E^k$ with $\{E^k\}$ being the frame in the
cotangent bundle dual to $\{E_k\}$, \ie  $E^k(E_i)=\delta^k_i$ with
$\delta_{i}^{j}$ being the Kronecker deltas
($\delta_{i}^{j}=0$ for $i\not=j$ and
($\delta_{i}^{j}=1$ for $i=j$);
$ \Hat{T}(u,\xi):= \nabla_u\xi - \nabla_\xi u - \nabla_{[u,\xi]} $
is the torsion operator in $L_n$;
\(
\Lied_\xi\Gamma(u,u)
:=
\Lied_\xi\Gamma_{ij}^{k} E^i\otimes E^j\otimes E_k
\)
with $\otimes$ being the tensor product sign;
 $\Lied_\xi u= (\Lied_\xi u^k) E_k = [\xi,u] = -\Lied_u \xi $.
	\end{Rem}

	\begin{Rem}	\label{Rem2.2}
Substituting in~\eref{1.4} the equality
\[
u^j \Lied\Gamma_{ji}^{k}
=
\Lied_\xi u_{|i}^{k} - [\Lied_\xi\bar{\nabla}_i] u^k
\]
we get get the following form of the generalized deviation equation:
	\begin{align}
\tag{\ref{1.4}a}	\label{1.4b}
\frac{\bar{D}^2\xi^k}{\od s^2}
& =
R_{ijl}^{k} u^iu^j\xi^l +
\xi_{|j}^{k} F^j +
u^j \frac{ \bar{D} (T_{jl}^{k}\xi^l) } {\od s} +
u^i \{ \Lied_\xi(u_{|i}^{k}) - (\Lied_\xi u^i)_{|i} \}
\\\intertext{or, in coordinate/frame independent notation,}
\tag{\ref{1.5}b}	\label{1.5b}
	\begin{split}
\frac{D^2\xi}{\od s^2}
& =
\Hat{R}(u,\xi)u + C_{1}^{1}(F\otimes D\xi) +
\frac{D(\Hat{T}(F,\xi))}{\od s} - \Hat{T}(F,\xi)
- [\nabla_u,\Lied_\xi] u .
	\end{split}
	\end{align}
	\end{Rem}


\section {Physical interpretation of the deviation equation}
\label{Sect2}

\indent\textbf{2.1.}
	In this section, when interpreting physically the deviation equation,
we shall restrict ourselves to infinitesimal deviation vector with components
$\xi^i$ defined below via~\eref{2.1a}, which means that we shall have in mind
the so\ndash called \emph{local} deviation equation. The cause for that is
in the non\ndash local problems arising when one tries to derive  non\ndash
local deviation equations, which are connected with the comparison of tensors
defined at non\ndash infinitesimally near points in $L_n$. These problems are
out of the range of this investigation and do not have a unique and global
solution at present. However, there are arguments that indicate the validity
in the nonlocal case of the physical interpretation of the deviation equation
presented below.

\textbf{2.2.}
	Consider two point-like particles 1 and 2 moving along the $C^2$
paths $x_1(\tau_1)$ and $x_2(\tau_2)$, respectively, in $L_n$ with parameters
$\tau_1$ and $\tau_2$, \ie $x_1(\tau_1)$ and $x_2(\tau_2)$ are their
trajectories. Assume the behaviour of the particles is ``observed'' by a
particle (observer) with $C^2$ trajectory $x_0(\tau)$, the basic trajectory,
with parameter $\tau$. It will be supposed that the mappings $x_0$, $x_1$ and
$x_2$ are injective and invertible on the sets of their ranges and that
$\tau\not=\const$, \ie $\od\tau\not=0$.

	Mathematically the ``observation'' process means existence of
mappings $\phi_1$ and $\phi_2$ such that, when the observer is at a point
$x_0(\tau)$ and the observed particles at $x_1(\tau_1)$ and $x_2(\tau_2)$, we
have
$x_1(\tau_1)=\phi_1(x_0(\tau))$ and
$x_2(\tau_2)=\phi_2(x_0(\tau))$.
This means that, when the observer ``knows'' its own position in $L_n$, he
can find the positions of the observed particles.

	Let us put:
	\begin{align}	\label{2.1}
u^\alpha &:= \frac{\od x_0^\alpha(\tau)}{\od \tau},\quad
u:=u^\alpha\pd_\alpha=u^i E_i, \quad
F:=\frac{Du}{\od\tau}=u^i\nabla_iu
\\ \tag{\ref{2.1}a}	\label{2.1a}
\xi^\alpha
&:= x_2^\alpha(\tau_2) - x_1^\alpha(\tau_1)
 = \phi_2^\alpha(x_0(\tau)) - \phi_1^\alpha(x_0(\tau))
  = \xi^\alpha(\tau)
\\ \tag{\ref{2.1}}	\label{2.1b}
\xi &:= \xi^\alpha\pd_\alpha=\xi^k E_k
\\ \tag{\ref{2.1}c}	\label{2.1c}
V &:= \frac{D\xi}{\od\tau} = u^k\nabla_k\xi,\quad
V^i=\frac{\bar{D}\xi^i}{\od\tau}
    = \frac{\od\xi^i}{\od\tau} +\Gamma_{jk}^{i}\xi^ju^k.
	\end{align}

	The contravariant vector $\xi$ will be called the \emph{deviation
vector}. Evidently, it describes the relative position of the particle 2 with
respect to the particle 1 (as seen from a point on the basic trajectory). The
vector $V$ is physically interpreted as a relative velocity of the particle 2
with respect to particle 1 as seen from the basic trajectory whose parameter
$\tau$ is interpreted as observer's ``proper'' time.

	The deviation equation in the form~\eref{1.5a} for the deviation
vector $\xi$ reads
	\begin{multline}	\label{2.2}
\frac{D^2\xi}{\od\tau^2}
= \frac{DV}{\od\tau}
=
\Hat{R}(u,\xi)u + C_1^1(F\otimes(D\xi)) + \frac{\Hat{T}(u,\xi))}{\od\tau}
\\
- \Hat{T}(F,\xi) + \Lied_\xi F
- \frac{D(\Lied_\xi u)}{\od\tau} - C_1^1((Du)\otimes(\Lied_\xi u)) .
	\end{multline}

	It is clear, with respect to a fixed observer, the vector
$\frac{DV}{\od\tau}$ is the relative acceleration between the particle 2
relative to the particle 1. Therefore the generalized deviation equation
gives the relative acceleration between the observed particles as a function
of the space $L_n$, \ie $\Hat{R}$, $\Hat{T}$, and $\Gamma_{jk}^{i}$, the
trajectory of the observer ($\tau$, $u$ and $F$), and the relative motion of
the observed particles ($\xi$ and $V$). It should be recalled, together
with~\eref{2.2}, in any particular case, one should consider also a suitable
additional condition(s).

	Notice, any type of restrictions on the trajectories considered can
be regarded as additional conditions to the deviation equation.


\section
[Examples of additional conditions defining deviation equations]
{Examples of additional conditions defining\\ deviation equations}
\label{Sect3}

\vspace{1.5ex}
\indent\textbf{3.1.}
\underline{The symmetries of $L_n$ as a source of additional conditions.}\\

	If the space $L_n$ admits some symmetries, this immediately entails
definite additional conditions in the generalized deviation
equation~\cite{Manov-	1978,Norden}. This is so because, when deriving the
deviation equation, we have used only quantities characterizing
$L_n$ without imposing on them any restrictions.

	To illustrate the above, we present bellow four examples. The
mathematical derivation of the corresponding additional conditions
in~\cite{Norden,Sinyukov,Schouten/physics,Yano/LieDerivatives} for Riemannian
space $V_n$ and for $L_n$ space it can be obtained in a way similar to the
one pointed in these references.

\textbf{(i)}
	If $L_n$ admits geodesic mappings, which map the geodesics of $L_n$
in geodesics if $L_n$, we have the equations
	\begin{gather}	\label{3.1}
\Lied_\xi \overset{\mathrm{p}}{\Gamma}{}_{jk}^{i} = 0
\\\intertext{where}
\tag{\ref{3.1}a}	\label{3.1a}
\overset{\mathrm{p}}{\Gamma}{}_{jk}^{i}
:=
\Gamma_{(jk)}^{i} -
\frac{2}{n+1} \bigl( \delta_{j}^{i} \Gamma_{(lk)}^{l} \bigr)_{(jk)}
	\end{gather}
are the Thomas projective parameters~\cite{Yano/LieDerivatives}
(symmetrization is performed  over the indices included in parentheses ,
$B_{(ij)} ;= \frac{1}{2}(B_{ij}+B_{ji})$);

\textbf{(ii)}
	If $L_n$ admits affine transformations, which preserve the the
parallelism of the tangent vectors, this leads to the conditions
	\begin{equation}	\label{3.2}
\Lied_\xi\Gamma_{jk}^{i} = 0 .
	\end{equation}

\textbf{(iii)}
	If in $L_n$ a metric with local components $g_{ij}(=g_{ji})$ is
given and $L_n$ admits symmetries, \ie mappings $L_n\to L_n$ preserving the
distances defined by the metric, than the conditions
	\begin{equation}	\label{3.3}
\Lied_\xi g_{ij} = 0
	\end{equation}
act as additional conditions to the generalized deviation equation.

\textbf{(iv)}
	If in $L_n$ a metric $g_{ij}$ is given and the space admits conformal
transformations, which preserve the angles between the tangent vectors, then
there exists a function $\Phi\colon L_n\to\field[R]$ such that
	\begin{equation}	\label{3.4}
\Lied_\xi g_{ij} = 2 \Phi g_{ij}.
	\end{equation}
If the metric tensor is non\ndash degenerate, $\det[g_{ij}]\not=0,\infty$,
the equation~\eref{3.4} can be rewritten as
	\begin{equation}
\tag{\ref{3.4}a}	\label{3.4a}
\Lied_\xi\bigl( |g|^{-1/n} g_{ij} \bigr) =0
	\end{equation}
with $g:=\det[g_{ij}]$ and $n$ being the dimension of the $L_n$-space.

\vspace{1.5ex}
\textbf{3.2.}
\underline{Particular realizations of $L_n$ as examples of additional
conditions.}\\

As additional conditions in the deviation equation can play role relations that
connect structures over $L_n$ or defining new objects over $L_n$ and their
possible connections with the already existing ones. Examples of this type are
presented in the following list:

\textbf{(i)}
	$L_n$-spaces without torsion,
    \begin{equation}    \label{3.5}
T_{jk}^{i} = 0 \qquad (\hat{T}=0) .
    \end{equation}

\textbf{(ii)}
	$p$-recurrent, $p\in\field[N]$, $L_n$-spaces,
    \begin{equation}    \label{3.6}
R_{jkl|i_1\dots i_p}^{i} = R_{jkl}^{i} \cdot A_{i_1\dots i_p}
    \end{equation}
with $A_{i_1\dots i_p}$ being the components of a tensor.

\textbf{(iii)}
	Affine (locally flat)  $L_n$-spaces: for every point $x_0\in L_n$,
there exists its neighborhood $U$ such that
    \begin{equation}    \label{3.7}
R_{jkl}^{i}|_U = 0 \qquad (\hat{R}|_U = 0)
    \end{equation}
or, equivalently, their is a frame $E_{j_0}$ such that
    \begin{equation}
\tag{\ref{3.7}a}	\label{3.7a}
\Gamma_{j_0k_0}^{i_0}|_U = 0 .
    \end{equation}

\textbf{(iv)}
	Equiaffine (generally with torsion) space: $L_n$-space with symmetric
Ricci tensor,
    \begin{equation}    \label{3.8}
R_{ij}=R_{ji}\qquad (R_{ij}:=R_{ijk}^{k}).
    \end{equation}

\textbf{(v)}
	A metrical $L_n$-space with metric tensor with components $g_{ij}$
and additional vector $w_i$ that satisfy the semi-metrical transport condition,
    \begin{equation}    \label{3.9}
g_{ij|k} = w_k \cdot g_{ij} .
    \end{equation}
In particular, of this kind are the Weil spaces in which the metric is
nondegenerate and the inverse metric tensor $g^{ij}$ is defined via
    \begin{equation}    \label{3.10}
g^{ik} g_{kj} = \delta_j^i .
    \end{equation}

\textbf{(vi)}
	Einstein spaces (generally with torsion): metrical $L_n$-spaces in
which, for some scalar function $f$,
    \begin{equation}    \label{3.11}
R_{ij} = f\cdot g_{ij} .
    \end{equation}

\textbf{(vii)}
	Riemannian manifolds (generally with torsion): equiaffine Weil spaces
with metrical transport:
    \begin{equation}    \label{3.12}
g_{ij|k} = 0 \qquad (g=\det[g_{ij}] \not=0,\infty;\ R_{ij}=R_{ji}) .
    \end{equation}

\textbf{(viii)}
	Conformal-Euclidean spaces (generally with torsion): Weil spaces in
which exists a tensor with components $P_{ij}$ such that
    \begin{align}    \label{3.13}
R_{jkl}^{i}
&= 2 (\delta_k^i P_{lj} - \delta_j^i P_{kl}
	- g_{jk} P_{lm} g^{mi} )_{[kl]}
\\
\tag{\ref{3.13}a}	\label{3.13a}
\bar{\nabla}_{[i} P_{j]k}
& = 0 .
    \end{align}

For $n\ge3$, from~\eref{3.13} follow~\eref{3.13a} and
    \begin{equation}    \label{3.13b}
P_{ij}
= \frac{1}{n-2}
\Bigl(
R_{ij} + \frac{2}{n} R_{[ij]} - \frac{1}{2(n-1)} (R_{kl} g^{lk}) g_{ij}
\Bigr) .
    \end{equation}

	Anyone of the conditions i--vii presented above, or a collection of them
(if they are compatible), together with the identity~\eref{1.4} define a
particular deviation equation in the corresponding spaces.

\vspace{1.5ex}
\textbf{3.3.}
\underline{Other additional conditions.}\\

	The study of concrete problems, connected with deviation equations, is
based on mathematical or physical reasonings and leads to appropriate additional
conditions. Examples of such conditions and their interpretation in coordinate
frames are presented in~\cite{Manov-1978}, where Riemannian spaces without
torsion are considered. Without repeating this reference, below we present some
examples of additional conditions and the corresponding to them deviation
equations~\eref{1.5} or~\eref{1.5a} in spaces $L_n$ with affine connection.

\textbf{(i)}
	 $F=0$ ($x(s)$ is a geodesic path in $L_n$):
\[
\frac{D^2\xi}{\od s^2}
=
\Hat{R}(u,\xi)u
+ \frac{D}{\od s}(\hat{T}(u,\xi) - \Lied_\xi u)
- C_1^2((Du)\otimes (\Lied_\xi u))  .
\]

\textbf{(ii)}
	 $\Lied_\xi u= 0$:
\[
\frac{D^2\xi}{\od s^2}
=
\Hat{R}(u,\xi)u
+ C_1^1(F\otimes(D\xi)) + \frac{D}{\od s} (T(u,\xi))
- \hat{T}(F,\xi) + \Lied_\xi F .
\]

\textbf{(iii)}
	 $u^i_{|k}=0$ and hence $Du=0$ and $F=0$:
\[
\frac{D^2\xi}{\od s^2}
=
\Hat{R}(u,\xi)u
+ \frac{D}{\od s} (T(u,\xi) \Lied_\xi u) .
\]

\textbf{(iv)}
	 $\Lied_\xi F=-F$:
\[
\frac{D^2\xi}{\od s^2}
=
\Hat{R}(u,\xi)u
+ C_1^1(F\otimes(D\xi)) - F + \frac{D}{\od s} (T(u,\xi)-\Lied_\xi u)
- \hat{T}(F,\xi) - C_1^2((Du)\otimes(\Lied_\xi u)) .
\]

\textbf{(v)}
	 $u=\xi$ and hence $\Lied_\xi u=[\xi,u]=0$:
\[
\frac{D^2\xi}{\od s^2}
=
\Lied_uF + C_1^1(F\otimes(Du)) + \hat{T}(u,F) .
\]

\textbf{(vi)}
	 $u^iu^j\Lied_\xi\Gamma_{ij}^{k} = - (F^k+\xi_{|j}^k F^j)$:
\[
\frac{D}{\od s}\Bigl( u + \frac{D\xi}{\od s} \Bigr)
=
\Hat{R}(u,\xi)u
+ \frac{D}{\od s} (T(u,\xi)) - \hat{T}(F,\xi) .
\]


\section {Examples of deviation equations in $L_n$}
\label{Sect4}

\indent\textbf{4.1.}
Let $x(s,q)$, $s,q\in \field[R]$, be a 2-parameter $C^2$  family of curves in a
space $L_n$; rigorously speaking~\cite{Mitskevich}, the parameter $q$ should
be replaced with $n-1$ independent parameters $q_1,\dots,q_{n-1}$ whose
employment does not change essentially the next considerations.

	Define the vectors $u$ and $v$ as tangent to respectively the $s$-curves
and $q$-curves,
    \begin{equation}    \label{4.1}
    \begin{split}
u:=u^iE_i=u^\alpha\pd_\alpha &\qquad u^\alpha=\frac{\pd x^\alpha(s,q)}{\pd s}
\\
v:=v^iE_i=v^\alpha\pd_\alpha &\qquad v^\alpha=\frac{\pd x^\alpha(s,q)}{\pd q} .
    \end{split}
    \end{equation}
Let $\delta q$ be an infinitesimal constant. The vector
    \begin{equation}    \label{4.2}
\xi := (\delta q) v
    \end{equation}
can be interpreted as a relative distance between the curves
 $x(s,q)$ and $x(s,q+\delta q)$. Now we want to find the relative acceleration
between these two cures relative for the former one. The solution of that
problem is given via the deviation equation. The  $C^2$ character of the family
$x$ serves as an additional condition as it leads to
    \begin{equation}    \label{4.3}
0
=
  \frac{\pd^2 x^\alpha(s,q) }{\pd q \pd s}
- \frac{\pd^2 x^\alpha(s,q) }{\pd s \pd q}
=
  \frac{\pd u^\alpha}{\pd q}
- \frac{\pd v^\alpha}{\pd s}
=
\Lied_v u^\alpha
=
- \Lied_u v^\alpha .
    \end{equation}
From~\eref{4.2} and~\eref{4.3}, we derive
    \begin{equation}    \label{4.4}
\Lied_\xi u = - \Lied_u \xi = 0
    \end{equation}
which can be rewritten as
    \begin{equation}    \label{4.5}
\frac{D\xi}{\od s} = C_1^1 (\xi\otimes(Du)) - \hat{T}(\xi,u).
    \end{equation}
Hence the deviation equation takes the form
    \begin{equation}    \label{4.6}
\frac{D^2\xi}{\od s^2}
=
\hat{R}(u,\xi)u + C_1^1(\xi\otimes (DF+T))
    \end{equation}
with
    \begin{equation}
\tag{\ref{4.6}a}	\label{4.6a}
T:= T_l^k E_k\otimes E^l
\qquad
T_l^k
=-u^n \nabla_n(T_{lj}^{k} u^j) + u^j T_{ji}^{k} (u_{|l}^i - T_{lm}^{i} u^m) .
    \end{equation}
The traditional method for obtaining this particular deviation equation
consists in the application of the operator $\frac{D}{\od}=u^k\nabla_k$
to~\eref{4.5} and a suitable changes in the result obtained.

\textbf{4.2.}
Consider two free point particles with $C^2$ trajectories $x_0(r)$ and
$x_1(r_1)$; the former will be interpreted as an observer and the latter one
as observed particle. The freedom means that the trajectories are geodesics
which implies the equation
    \begin{gather}    \label{4.7}
F^i = u_{|i}^k u^k
=
A_\alpha^i
\Bigl(
\frac{\od^2x_0^\alpha(r) }{\od r^2 }
+ \Gamma_{\beta\gamma}^{\alpha}
\cdot \frac{\od x_0^\beta}{\od r} \frac{\od x_0^\gamma}{\od r}
\Bigr)
=
0
\\		    \label{4.8}
\Lied_\xi u = 0 .
    \end{gather}

	The condition~\eref{4.7} defines the geodesic $x_0(r)$,
while~\eref{4.8}, which is equivalent to~\eref{4.4}, characterizes the
existence of the mapping $\phi_2$ (see section~\ref{Sect2}); the
condition~\eref{4.7} can also be considered as obtained via the method in
the previous example (where the curves $x_0(r)$ and $x_1(r_1)$, belonging to the
family $x(s)$ for $s=r,r_1$, are studied (for $q=\const=q_0$, the curves
$x(s=r,q_0)$ are geodesics).

	Now the deviation equation takes the following form (see point~3.3,
case~2 above):
    \begin{equation}    \label{4.9}
\frac{DV}{\od r} = \hat{R}(u,\xi)u + C_1^1(\xi\otimes T) ,
    \end{equation}
where
    \begin{equation}    \label{4.10}
V:=\frac{D\xi}{\od r} \quad
T:= T_l^k E_k^l \quad
T_l^k = -T_{lj|n}^{k}u^ju^n + u^jT_{ji}^{k}(u_{|l}^i - T_{lm}^{i}u^m) .
    \end{equation}

	From~\eref{4.9}, we conclude that the relative acceleration between two
free particles (``tidal'' acceleration) has two sources: the curvature tensor
and the torsion tensor. For that reason, in a $L_n$-space is possible a
situation which cannot arise in the $V_p$-spaces (without torsion): the torsion
can `compensate' the curvature so that the tidal acceleration will vanish,
$\frac{D V }{\od r}=0$, along the whole base trajectory, on some its parts, or
at some single point on it. Consequently, if one can experimentally prove that
in a non-flat spacetime region is realized one of these possibilities, this
possible fact will mean that in the region mentioned the torsion does not
vanish (as $L_n$ is used as a spacetime model).

\textbf{4.3.}
	Let an observer and observed particle have arbitrary $C^2$ trajectories
$x_0(r)$ and $x_1(r_1)$, respectively. Let us define
    \begin{equation}    \label{4.11}
    \begin{split}
& w:=\frac{\od r_1}{\od r}\ (\not=0,\infty) \quad
u^\alpha:=\frac{\od x_0^\alpha(r)}{\od r} \quad
u:=u^a\pd_\alpha=u^kE_k
\\
& v^\alpha:=\frac{\od x_1^\alpha(r_1)}{\od r_1} \quad
v:=v^a\pd_\alpha=v^kE_k .
    \end{split}
    \end{equation}
As $x_0(r)$ and $x_1(r_1)$ are, by definition, bijective and invertible
functions of $r$ and $r_1$, respectively, we can write
    \begin{equation}    \label{4.12}
v=\frac{1}{w} (u+V) \qquad
V:=\frac{D\xi}{\od r} .
    \end{equation}

	The mapping $x_0(r)\mapsto x_1(r_1)=x_0(r)+\xi(r)$ deforms the vector
$F^k:=u_{|l}^ku^l$ and the connection $\Gamma_{jk}^{i}$ in respectively
    \begin{align}    \label{4.13}
\lindex[F]{}{\prime}^i
& = F^i + \varepsilon\Lied_\zeta F^i
= F^i + \Lied_\xi F^i
\\		    \label{4.14}
\lindex[\Gamma]{}{\prime}_{jk}^{i}
& = \Gamma_{jk}^{i} + \varepsilon\Lied_\zeta \Gamma_{jk}^{i}
= \Gamma_{jk}^{i} + \Lied_\xi \Gamma_{jk}^{i}
    \end{align}
A natural additional condition for the deviation equation, in the particular
case, can be obtained by requiring $\lindex[F]{}{\prime}^i$ to describe the
change of $v^i$ on the observed trajectory, \viz
    \begin{equation}    \label{4.15}
\lindex[F]{}{\prime}^i
=v^j \lindex[\nabla]{}{\prime}_j v^i
= v^j( v_{,j}^i + \lindex[\Gamma]{}{\prime}_{jk}^{i} v^k) .
    \end{equation}
This means that the mapping $x_0(r)\mapsto x_1(r_1)$ draggs/deforms all
structures on $x_0(r)$ and the links between them to similar ones but defined
on the observed trajectory $x_1(r_1)$.

	One can prove that a necessary and sufficient condition for the
validity of~\eref{4.15} is
    \begin{multline}    \label{4.16}
\Lied_\xi F^i
=
-F^i + \frac{1}{w^2}
\Bigl\{
F^i - (u^i+ V^i) \Bigl( V^j(\ln w)_{,j} \frac{D \ln w}{\od r} \Bigr)
\\+ V^j (u^j+V^i)_{,j} + u^j V_{|j}^i
+ (u^j+V^j)(u^k+V^k) \Lied_\xi\Gamma_{jk}^{i}
\Bigr\} .
    \end{multline}

	In a similarly way one can derive other additional
condition~\cite{Manov-1978} if, instead of $v^\alpha:=\frac{\od
x_1^\alpha(r_1)}{\od r_1}$, one defines $v^\alpha=\frac{1}{w}u^\alpha$ and
requires the fulfillment of~\eref{4.15}. In this case, the transformation
$x_0(r)\mapsto x_1(r_1)$ maps/deforms the basic trajectory in such a way that
the tangent vectors to basic and deformed (observed) curves are proportional.

\vspace{1ex}
	The authors are thankful to N.~A.~Charnikov for discussions.


\renewcommand{\thesection}{A}
\section {Appendix: Lie derivatives of the coefficients
\protect{$\Gamma_{jk}^{i}$} of\\
 a connection in \protect{$L_n$} in arbitrary frame}

	 \label{Appendix}

	The Lie derivatives of the coefficients $\Gamma_{jk}^{i}$ of an affine
connection in $L_n$ are defined in a coordinated frame, e.g.,
in~\cite{Yano/LieDerivatives}. Below only the case of an arbitrary frame
$\{E_i\}$ will be considered.
	Recall, under a frame change $\{E_i\}\mapsto\{E_{i'}:=A_{i'}^jE_j\}$,
with $A:=[A_{i'}^i]$ being a nondegenerate matrix-valued function, the
$\Gamma$'s transform according to
    \begin{equation}    \label{A.1}
\Gamma_{jk}^{i} \mapsto \Gamma_{j'k'}^{i'}
=
A_{i}^{i'} A_{j'}^{j} A_{k'}^{k} \Gamma_{jk}^{i} + A_{i}^{i'} A_{j'k'}^{i} ,
    \end{equation}
where $[A_{i}^{i'}]:=A^{-1}$.

	Let us consider in $L_n$ an infinitesimal point transformation
    \begin{equation}    \label{A.2}
x\mapsto \bar{x}
    \end{equation}
with coordinate representation
    \begin{equation}    \label{A.3}
\bar{x}^\alpha = x^\alpha + \varepsilon \xi^\alpha(x^\beta) .
    \end{equation}
Here $\varepsilon$ is an infinitesimal parameter and $\xi^\alpha$ are the
components of  $C^2$ vector field. The frame
$\{E_{\bar{i}}=A_{\bar{i}}^{i}E_i\}$ with
    \begin{equation}    \label{A.4}
A_{\bar{i}}^{j} = \delta_{\bar{i}}^{k} ( \delta_k^j - \varepsilon \Sigma_k^j )
\qquad
\Sigma_k^j := \xi_{,k}^j + C_{kl}^{j}\xi^l
\quad
\xi^l := A_\alpha^l \xi^\alpha
    \end{equation}
is called dragged by the point transformation~\eref{A.2}
(cf.~\cite{Schouten/physics}). The so-defined quantities satisfy the equations
    \begin{equation}    \label{A.5}
\Lied_\xi E_i = - \Sigma_i^kE_k
\quad
\delta_i^{\bar{j}} E_{\bar{j}}
= E_i +\varepsilon\Lied_\xi E_i
= \delta_i^{\bar{j}} A_{\bar{j}}^{k} E_k .
    \end{equation}
  Equation~\eref{A.4} implies
    \begin{equation}    \label{A.6}
A_j^{\bar{i}} = \delta_k^{\bar{i}} (\delta_j^k + \varepsilon\Sigma_j^k),
    \end{equation}
where here and henceforth the terms of order $\varepsilon^2$ and more are
neglected.

	By definition, the dragged (deformed) by~\eref{A.2} coefficients
$\lindex[\Gamma]{}{\prime}_{jk}^{i}$ are (cf.~\cite{Yano/LieDerivatives})
    \begin{equation}    \label{A.7}
\lindex[\Gamma]{}{\prime}_{jk}^{i}(x)
=
A_{\bar{i}}^{i} A^{\bar{j}}_{j}A^{\bar{k}}_{k}
\delta^{\bar{i}}_{p} A_{\bar{j}}^{q}A_{\bar{k}}^{r}
\Gamma_{qr}^{p}(x)
+
A_{\bar{i}}^{i}(x) A^{\bar{i}}_{j,k}(x) .
    \end{equation}
Inserting~\eref{A.4} and~\eref{A.6}  into the last expression and using that
    \begin{equation}    \label{A.9}
\Gamma_{qr}^{p}(\bar{x})
=
\Gamma_{qr}^{p}(x) + \varepsilon\xi^m \Gamma_{qr,m}^{p}(x),
    \end{equation}
we get
    \begin{equation}    \label{A.10}
\lindex[\Gamma]{}{\prime}_{jk}^{i}
=
\Gamma_{jk}^{i}
+ \varepsilon \bigl(
- \Gamma_{jk}^{n} \Sigma_n^i
+ \Gamma_{nk}^{i} \Sigma_j^n
+ \Gamma_{jn}^{i} \Sigma_k^n
+ \Sigma_{j,k}^i
+  \Gamma_{jk,n}^{i} \xi^n
\bigr) .
    \end{equation}
The Lie derivative of $\Gamma_{jk}^{i}$ with respect to $\xi$
is~\cite{Yano/LieDerivatives}
    \begin{equation}    \label{A.11}
\Lied_\xi \Gamma_{jk}^{i}
:=
\lim_{\varepsilon\to0}
\frac{ \lindex[\Gamma]{}{\prime}_{jk}^{i} - \Gamma_{jk}^{i} }
     {\varepsilon} ,
    \end{equation}
from where one immediately finds that (see~\eref{A.10})
    \begin{equation}    \label{A.12}
\Lied_\xi \Gamma_{jk}^{i}
=
 - \Gamma_{jk}^{n} \Sigma_n^i
+ \Gamma_{nk}^{i} \Sigma_n^j
+ \Gamma_{jn}^{i} \Sigma_n^k
+ \Sigma_{j,k}^i
+  \Gamma_{jk,n}^{i} \xi^n .
    \end{equation}

	If the partial derivatives in~\eref{A.12} are replace via covariant
ones ($\xi_{,k}^i=\xi_{|k}^i-\Gamma_{jk}^{i}\xi^j$), one gets the
equality~\eref{1.1} after appropriate calculations.


\addcontentsline{toc}{section}{References}
\bibliography{bozhopub,bozhoref}

\begin{thebibliography}{1}

\bibitem{Manov-1978}
Sawa~S. Manov.
\newblock The deviation equations and {Lie} derivatives in {Riemannian} spaces.
\newblock JINR Communication P2-12026, Dubna, 1978.
\newblock In Russian.

\bibitem{MTW}
C.~W. Misner, K.~S. Thorne, and J.~A. Wheeler.
\newblock {\em Gravitation}.
\newblock W.~H. Freeman and Company, San Francisco, 1973.

\bibitem{Mitskevich}
N.~V. Mitskevich.
\newblock {\em Physical fields in general theory of relativity}.
\newblock Nauka, Moscow, 1969.
\newblock (In Russian).

\bibitem{Yano/LieDerivatives}
K.~Yano.
\newblock {\em The theory of Lie derivatives and its applications}.
\newblock North-Holland Publ.\ Co., Amsterdam, 1957.

\bibitem{Norden}
A.~P. Norden.
\newblock {\em Spaces with affine connection}.
\newblock Nauka, Moscow, second edition, 1976.
\newblock (In Russian).

\bibitem{Sinyukov}
N.~S. Sinyukov.
\newblock {\em Geodesic mappings of Riemannian manifolds}.
\newblock Nauka, Moscow, 1979.
\newblock In Russian.

\bibitem{Schouten/physics}
J.~A. Schouten.
\newblock {\em Tensor analysis for physicists}.
\newblock Clarendon Press, Oxford, 1951.

\end{thebibliography}
\bibliographystyle{unsrt}
\addcontentsline{toc}{subsubsection}{This article ends at page}

\end{document}